\documentclass[12pt]{amsart}

\usepackage{geometry}

\usepackage{color}

\numberwithin{equation}{section} 
\numberwithin{figure}{section} 
\theoremstyle{plain}
\newtheorem{thm}{Theorem}[section]
\theoremstyle{plain}
\newtheorem{lem}[thm]{Lemma}
\theoremstyle{plain}
\newtheorem{prop}[thm]{Proposition}
\theoremstyle{plain}
\newtheorem{cor}[thm]{Corollary}
\theoremstyle{remark}
\newtheorem{rem}[thm]{Remark}
\theoremstyle{remark}

\theoremstyle{remark}
\newtheorem{example}[thm]{Example}
\theoremstyle{plain}
\newtheorem*{prop*}{Proposition}
\newtheorem*{thm*}{Theorem}

\theoremstyle{remark}
\newtheorem*{rem*}{Remark}

\theoremstyle{plain}
\newtheorem{question}{Question}

\theoremstyle{remark}

\theoremstyle{definition}

\theoremstyle{definition}
\newtheorem{defi}[thm]{Definition}

\usepackage{latexsym,amsfonts, xypic,amssymb}
\xyoption{all}

\newcommand{\C}{\mathbb{C}}

\newcommand{\Z}{\mathbb{Z}}

\newcommand{\N}{\mathbb{N}}

\newcommand{\der}{\partial}

\newcommand{\Gr}{\textup{Gr}}
\newcommand{\spec}{\textup{Spec}}
\newcommand{\LND}{\textup{LND}}
\newcommand{\ML}{\textup{ML}}

\newcommand{\Ker}{\textup{Ker}}

\begin{document}

\title{Classification(s) of Danielewski hypersurfaces}

\author{P.-M. Poloni}
\address{Mathematisches Institut
Universit\"at Basel, Rheinsprung 21, CH-4051 Basel, Switzerland}

\email{pierre-marie.poloni@unibas.ch }

\begin{abstract}
The Danielewski hypersurfaces are the hypersurfaces $X_{Q,n}$ in
$\C^3$ defined by an equation of the form $x^ny=Q(x,z)$ where
$n\geq1$ and $Q(x,z)$ is a polynomial such that $Q(0,z)$ is of
degree at least two. They were studied by many authors during the
last twenty years. In the present article, we give their
classification as algebraic varieties. We also give their
classification up to automorphism of the ambient space. As a
corollary, we obtain that every Danielewski hypersurface $X_{Q,n}$
with $n\geq2$ admits at least two non-equivalent embeddings into
$\C^3$.
\end{abstract}

\maketitle

\section{Introduction}

The story of Danielewski hypersurfaces goes back to 1989, when
Danielewski \cite{Danielewski} showed that, if $W_n$ denotes the
hypersurface in $\C^3$ defined by the equation $x^ny-z(z-1)=0$,
then $W_n\times\C$  and $W_m\times\C$  are isomorphic algebraic
varieties for all $n,m\geq1$, whereas the surfaces $W_1$ and $W_2$
are not isomorphic. He discovered  the first counterexamples to
the Cancellation Problem. Then, Fieseler \cite{Fieseler} proved
that $W_n$ and $W_m$ are not isomorphic if $n\neq m$.

Since these results appeared, complex algebraic surfaces defined
by equations of the form $x^ny-Q(x,z)=0$ (now called
\emph{Danielewski hypersurfaces}) have been studied by many
different authors (see  \cite{Wilkens}, \cite{ML01},
\cite{Daigle}, \cite{Crachiola}, \cite{Freudenburg-MJ},
\cite{MJ-P}), leading to new interesting examples as byproducts.
Let us mention two of them.

In their work on embeddings of Danielewski hypersurfaces given by
$x^ny=p(z)$, Freudenburg and Moser-Jauslin \cite{Freudenburg-MJ}
discovered  an example of two smooth algebraic surfaces which are
algebraically non-isomorphic but holomorphically isomorphic.

More recently, the study of Danielewski hypersurfaces of equations
$x^2y-z^2-xq(z)=0$ produced the first counterexamples to the
Stable Equivalence Problem \cite{MJ-P}; that is two polynomials of
$\C[X_1,X_2,X_3]$ which are not equivalent (i.e. such that there
exist no algebraic automorphism of $\C[X_1,X_2,X_3]$ which maps
one to the other one) but, when considered as polynomials of
$\C[X_1,X_2,X_3,X_4]$, become equivalent.

The purpose of the present paper is  to classify all Danielewski
hypersurfaces, both as algebraic varieties, and also as hypersurfaces
in $\C^3$. More precisely, we will give necessary and sufficient
conditions for isomorphism of two Danielewski hypersurfaces; and,
on the other hand, we will give  necessary and sufficient
conditions for equivalence of two isomorphic Danielewski
hypersurfaces. Recall that two isomorphic hypersurfaces
$H_1,H_2\subset\C^n$ are said to be \emph{equivalent} if there
exists an algebraic automorphism $\Phi$ of $\C^n$ which maps one
to the other one, i.e. such that $\Phi(H_1)=H_2$.

We know indeed that isomorphic classes and equivalence classes are
distinct for Danielewski hypersurfaces. This was first observed by  Freudenburg and
Moser-Jauslin, who  showed in \cite{Freudenburg-MJ} that the
Danielewski hypersurfaces defined respectively by the equations
$f=x^2y-(1+x)(z^2-1)=0$ and $g=x^2y-z^2+1=0$ are isomorphic but
non equivalent. (One way to see that they are not equivalent is to
remark that the level surfaces $f^{-1}(c)$ are smooth for every
constant $c\in\C$, whereas the surface $g^{-1}(1)$ is singular
along the line $\{x=z=0\}$.)

Several papers already contain the classification, up to isomorphism,
of Danielewski hypersurfaces of a certain form. Makar-Limanov
proved in \cite{ML01} that two Danielewski hypersurfaces of
equations $x^{n_1}y-p_1(z)=0$ and $x^{n_2}y-p_2(z)=0$ with
$n_1,n_2\geq 2$ and $p_1,p_2\in\C[z]$ are isomorphic if and only
if they are equivalent via an affine automorphism of the form
$(x,y,z)\mapsto(ax,by,cz+d)$ with $a,b,c\in\C^*$ and $d\in\C$.
Then, Daigle  generalized in \cite{Daigle} this result to the case
$n_1,n_2\geq1$. Next, Wilkens  has given in \cite{Wilkens} the
classification of Danielewski hypersurfaces of equations
$x^{n}y-z^2-h(x)z=0$ with $n\geq2$ and $h(x)\in\C[x]$.

Finally, Dubouloz and the author showed in \cite{Dubouloz-Poloni}
that every Danielewski hypersurface $X_{Q,n}$ of equation
$x^ny=Q(x,z)$, where $Q(x,z)$ is such that $Q(0,z)$ has  simple
roots, is isomorphic to a one defined by an equation of the form
$x^ny=\prod_{i=1}^d(z-\sigma_i(x))$, where
$\{\sigma_1(x),\ldots\sigma_d(x)\}$ is a collection of polynomials
in $\C[x]$ so that $\sigma_i(0)\neq\sigma_j(0)$ if $i\neq j$. In
the same paper, we classified these last ones and  called them
\emph{standard forms}. This effectively classifies, up to
isomorphism, all Danielewski hypersurfaces of equations
$x^ny=Q(x,z)$, where  $Q(0,z)$ has  simple roots.

In the present paper, we  generalize the notion of Danielewski
hypersurface in standard form and we prove that every Danielewski
hypersurface is isomorphic to a one in standard form (which can be
found by an algorithmic procedure). Then,  we are able to classify
all Danielewski hypersurfaces. The terminology \emph{standard
form} is relevant since every isomorphism between two Danielewski
hypersurfaces in standard form -- and every automorphism of such a
Danielewski hypersurface -- extend to a triangular automorphism of
$\C^3$.

We also give a criterion (Theorem \ref{thm:classification-plgts}) to distinguish
isomorphic but not equivalent Danielewski hypersurfaces.

As a corollary, we obtain that every Danielewski hypersurface
defined by an equation of the form $x^ny-Q(x,z)=0$ with $n\geq2$
admits at least two non equivalent embeddings into $\C^3$.

Most of these results are based on a precise picture of the sets
of locally nilpotent derivations of coordinate rings of
Danielewski hypersurfaces, obtained using techniques which were
mainly developed by Makar-Limanov in \cite{ML01}.

The paper is organized as follows. Section 1 is the introduction.
In section 2, we  fix some notations and definitions. In section
3, we study the locally nilpotent derivations on the Danielewski
hypersurfaces in order to get information on what an isomorphism
between two Danielewski hypersurfaces looks like.  Section 4 is
devoted to the classification of Danielewski hypersurfaces up to
equivalence, whereas sections 5 and 6 contain their classification
up to isomorphism and the study of the Danielewski hypersurfaces
in standard form.

\section{Definitions and notations}

In this paper, our base field is $\C$, the field of complex numbers. If $n\geq1$, then $\C^{[n]}$ will denote a polynomial ring in $n$ variables over $\C$.

\begin{defi}
Two hypersurfaces $X_1$ and $X_2$ of  $\C^n$ are said to be \emph{equivalent} if there exists a (polynomial) automorphism $\Phi$ of $\C^n$ such that $\Phi(X_1)=X_2$.
\end{defi}

This notion is related to the notion of \emph{equivalent embeddings} in the following sense. If $X_1$ and $X_2$ are two isomorphic hypersurfaces of $\C^n$ which are not equivalent, then $X_1$ admits two non-equivalent embeddings into $\C^n$. More precisely,  let $\varphi:X_1\to X_2$ be an isomorphism and denote $i_1:X_1\to\C^n$ and $i_2:X_2\to\C^n$ the inclusion maps. Then, $i_1$ and $i_2\circ\varphi$ are two non-equivalent embeddings of $X_1$ into $\C^n$, since $\varphi$ does not extend to an automorphism of $\C^n$.

\begin{defi}\label{defiHypDan}
A \emph{Danielewski hypersurface} is a hypersurface
$X_{Q,n}\subset\C^3$ defined by an equation of the form
$x^ny-Q(x,z)=0$, where $n\in\N$ and $Q(x,z)\in\C[x,z]$  is such
that $\deg\left(Q\left(0,z\right)\right)\geq2$.

We will denote by $S_{Q,n}$ the coordinate ring of a Danielewski
hypersurface $X_{Q,n}$, i.e.
$S_{Q,n}=\C[X_{Q,n}]=\C[x,y,z]/(x^ny-Q(x,z))$.
\end{defi}

It can be easily seen that every Danielewski hypersurface is
equivalent to a one of the form $X_{Q,n}$ with $\deg_xQ(x,z)<n$.

\begin{lem}\label{lem:mod_x^n}
Let $X_{Q,n}$ be a Danielewski hypersurface and $R(x,z)\in\C[x,z]$
be a polynomial. Then $X_{Q,n}$ is equivalent to the Danielewski
hypersurface of equation $x^ny-Q(x,z)-x^nR(x,z)$.
\end{lem}

\begin{proof}
It suffices to consider the  triangular automorphism of $\C^3$
defined by $(x,y,z)\mapsto(x,y+R(x,z),z)$.
\end{proof}

\section{Using locally nilpotent derivations}

One important property of Danielewski hypersurfaces is that they
admit nontrivial actions of the additive group $\C_+$. For
instance, we can define a $\C_+$-action $\delta_{Q,n}:\C\times
X_{Q,n}\to X_{Q,n}$ on a hypersurface $X_{Q,n}$ by posing
$$\delta_{Q,n}\left(t,\left(x,y,z\right)\right)=\left(x,y+x^{-n}\left(Q\left(x,z+tx^n\right)-Q\left(x,z\right)\right),z+tx^n\right).$$

Since a $\C_+$-actions on an affine complex surface $S$ induces a
$\C$-fibration  over an affine curve, affine complex surfaces with
$\C_+$-actions split into two cases. Either there is  only one
$\C$-fibration on $S$ up to an isomorphism of the base, or there
exists a second one. In other words, either the surface has a
Makar-Limanov invariant of transcendence degree one, or its
Makar-Limanov invariant is trivial. Recall that algebraic
$\C_+$-actions on an affine variety $\spec A$ correspond to
locally nilpotent derivations on the $\C$-algebra $A$ (for
example, the action $\delta_{Q,n}$ on a surface $X_{Q,n}$
corresponds to the locally nilpotent derivation
$\Delta_{Q,n}=x^n\frac{\der}{\der z}+\frac{\der
Q\left(x,z\right)}{\der z}\frac{\der}{\der y}$ on the coordinate
ring $\C[X_{Q,n}]$), and that the Makar-Limanov invariant $\ML(A)$
of an algebra $A$ is defined as the intersection of all kernels of
locally nilpotent derivations of $A$. Equivalently, $\ML(A)$ is
the intersection of all invariant rings of algebraic
$\C_+$-actions on $\spec(A)$.

Applying Makar-Limanov's techniques, one can obtain the first
important result concerning Danielewski hypersurfaces: the
Makar-Limanov invariant of a Danielewski hypersurface $X_{Q,n}$ is
non-trivial if $n\geq2$.

\begin{thm}\label{thm:ML}
Let $X_{Q,n}$ be a Danielewski hypersurface.  Then
$\ML(X_{Q,n})=\C$ if $n=1$ and $\ML(X_{Q,n})=\C[x]$ if $n\geq2$.
\end{thm}

\begin{proof}
Let $X_{Q,n}$ be a Danielewski hypersurface.

If $n=1$, the result is easy. Indeed, we can suppose that
$Q(x,z)=p(z)\in\C[z]$ (see Lemma \ref{lem:mod_x^n}). Then, it
suffices to consider the following locally nilpotent derivations
on the coordinate ring $S_{p,1}=\C[x,y,z]/(xy-p(z))$.
$$\delta_1=x\frac{\der}{\der z}+p'(z)\frac{\der}{\der y}\quad\text{ and }\quad \delta_2=y\frac{\der}{\der z}+p'(z)\frac{\der}{\der x}.$$
Since $\Ker(\delta_1)\cap\Ker(\delta_2)=\C$, this shows that the
Makar-Limanov invariant of every Danielewski hypersurface
$X_{Q,1}$ is trivial.

Suppose now that $n\geq2$ and let $\delta$ be a non-zero locally
nilpotent derivation on the coordinate ring $S_{Q,n}=\C[X_{Q,n}]$.
Without loss of generality, we can suppose that the leading term
of $Q(0,z)$ is $z^d$ with $d\geq2$.

Then, the proof given by Makar-Limanov in \cite{ML01} for hypersurfaces of equation
$x^ny=p(z)$  still holds.  This proof goes as follows.

The main idea is to consider $S_{Q,n}$ as a subalgebra of
$\C[x,x^{-1},z]$ with $y=x^{-n}Q(x,z)$ and to choose a
$\Z$-filtration on $S_{Q,n}$  such that the corresponding graded
algebra $\text{Gr}(S_{Q,n})$ is isomorphic to the subalgebra
$\C[x,x^{-n}z^d,z]$.

Recall that Makar-Limanov has proved that any non-zero locally
nilpotent derivation $D$ on an algebra $A$ with a $\Z$-filtration induces a
non-zero  locally nilpotent derivation $gr(D)$ on  the
graded algebra $\Gr(A)$. In our case,  he proved in \cite{ML01}
that $\Ker(gr(\delta))=\C[x]$.

Recall also that we can define a degree function associated to a
locally nilpotent $D\in\LND(A)$ by posing $\deg(0)=-\infty$ and
$\deg_D(a):=\max\{n\in\N\mid D^n(a)\neq0\}$ if $a\in
A\setminus\{0\}$. Moreover, if $A$ has a $\Z$-filtration, then
$\deg_{gr(D)}(gr(a))\leq\deg_D(a)$ for any element $a\in A$. (Here
$gr:A\to\Gr(A)$ denote the natural function from $A$ to $\Gr(A)$.)

This implies that $\Ker(\delta)=\C[x]$. Indeed, if $p\in
S_{Q,n}\setminus\C[x]$, then we can choose a filtration on
$S_{Q,n}$ such that $gr(p)$ belongs to
$\Gr(S_{Q,n})\setminus\C[x]$. Then $\Ker(\delta)\subset\C[x]$
follows from the inequalities
$1\leq\deg_{gr(\delta)}(gr(p))\leq\deg_\delta(p)$. Since
$\Ker(\delta)$ is of transcendence degree one and is algebraically
closed, we obtain that $x\in\Ker(\delta)$.

Thus, $\Ker(\delta)=\C[x]$ for any non-zero locally nilpotent
derivation on $S_{Q,n}$. In particular, $\ML(X_{Q,n})=\C[x]$.
\end{proof}

Using this result, we can obtain a precise picture of the set of locally
nilpotent derivations on rings $S_{Q,n}$ when $n\geq2$.

\begin{thm}\label{thm_LND}
Let $X_{Q,n}$ be a Danielewski hypersurface with  $n\geq2$ and let
$S_{Q,n}$ denote its coordinate ring. Then
$$\LND\left(S_{Q,n}\right)=\left\{h\left(x\right)\left(x^n\frac{\der}{\der
z}+\frac{\der Q\left(x,z\right)}{\der z}\frac{\der}{\der
y}\right), \text{ where } h(x)\in\C[x]\right\}.$$

In particular,  $\Ker(\delta)=\C[x]$ and
$\Ker(\delta^2)=\C[x]z+\C[x]$ for   every non-zero locally
nilpotent derivation $\delta\in\LND(S_{Q,n})$.
\end{thm}

\begin{proof}
Let $\delta$ be a non-zero locally nilpotent derivation on
$S_{Q,n}$ with $n\geq2$.  In the proof of Theorem \ref{thm:ML}, we
showed that  $\Ker(\delta)=\C[x]$. Then, due to Lemma 1.1 in
\cite{DaigleLND}, there exist polynomials $a(x),b(x)\in\C[x]$ such
that $a(x)\delta=b(x)\Delta_{Q,n}$ where
$\Delta_{Q,n}\in\LND(S_{Q,n})$ is the derivation defined by
$\Delta_{Q,n}=x^n\frac{\der}{\der z}+\frac{\der
Q\left(x,z\right)}{\der z}\frac{\der}{\der y}$.

The theorem will follow easily. First note that
$a(x)\delta(z)=x^nb(x)$. Also $a(x)$ divides $x^nb(x)$.
Therefore, in order to prove that $a(x)$ divides $b(x)$, it is
enough to show that $a(0)=0$ implies $b(0)=0$. This holds since
$\delta(y)\in\C[x,z]$ and
$a(0)\delta(y)(0,z)=b(0)(Q(0,z))'(z)$.

The theorem is proved.
\end{proof}

This theorem gives us a very powerful tool for classifying
Danielewski hypersurfaces. Indeed, note that
 an isomorphism $\varphi:A\to B$,  between two algebras $A$ and
$B$, conjugates the sets $\LND(A)$ and $\LND(B)$ of locally nilpotent derivations
on $A$ and
$B$, i.e. $\LND(A)=\varphi^{-1}\LND(B)\varphi$ if $\varphi:A\to B$
is an isomorphism. In turn, we obtain the following result.

\begin{cor}\label{cor:LND}\hfill
\begin{enumerate}
\item Let $\varphi:X_{Q_1,n_1}\to X_{Q_2,n_2}$ be an isomorphism
between two Danielewski hypersurfaces with $n_1,n_2\geq2$. Then,
there exist two constants $a,\alpha\in\C^*$ and a polynomial
$\beta(x)\in\C[x]$ such that $\varphi^*(x)=ax$ and
$\varphi^*(z)=\alpha z+\beta(x)$. \item If $X_{Q_1,n_1}$ and
$X_{Q_2,n_2}$ are two isomorphic Danielewski hypersurfaces, then
$n_1=n_2$ and $\deg(Q_1(0,z))=\deg(Q_2(0,z))$. \item Suppose that
$X_{Q_1,n}$ and $X_{Q_2,n}$ are two equivalent Danielewski
hypersurfaces with $n\geq2$, and let $\Phi:\C^3\to\C^3$ be an
algebraic automorphism such that $\Phi(X_{Q_1,n})=X_{Q_2,n}$.
Then, there exist constants $a,\alpha\in\C^*$, $\beta\in\C$ and a
polynomial $B\in\C^{[2]}$ such that $\Phi^*(x)=ax$ and
$\Phi^*(z)=\alpha z+\beta+xB(x,x^ny-Q_1(x,z))$.
\end{enumerate}
\end{cor}

\begin{proof}
For (1) and (2), we  follow the ideas of a proof given by Makar-Limanov in
\cite{ML01}.

Let $\varphi:X_{Q_1,n_1}\to X_{Q_2,n_2}$ be an isomorphism between
two Danielewski hypersurfaces with $n_1,n_2\geq2$. Let
$x_i,y_i,z_i$ denote the images of $x,y,z$ in the coordinate ring
$S_i=S_{Q_i,n_i}=\C[X_{Q_i,n_i}]$.

If $\delta\in\LND(S_1)$, then
$(\varphi^*)^{-1}\circ\delta\circ\varphi^*\in\LND(S_2)$. Thus,
Theorem \ref{thm_LND} implies  $\delta^2(z_1)=0$ and
$\delta^2\left(\varphi^*(z_2)\right)=0$ for any locally nilpotent
derivation $\delta\in\LND(S_1)$. Therefore,
$\varphi^*(z_2)=\alpha(x_1)z_1+\beta(x_1)$ for some polynomials
$\alpha$ and $\beta$. Since $\varphi$ is invertible, $\alpha$ must
be a nonzero constant $\alpha\in\C^*$.

On the other hand, $\varphi^*$ induces an isomorphism
$\ML(S_2)=\C[x_2]\to\ML(S_1)=\C[x_1]$. Consequently,
$\varphi^*(x_2)=ax_1+b$ for some constants $a\in\C^*$ and
$b\in\C$.

In order to prove $b=0$, consider the locally nilpotent derivation
$\delta_0\in\LND(S_2)$ defined by
$\delta_0=(\varphi^*)^{-1}\circ\left(x_1^n\frac{\der}{\der
z_1}+\frac{\der Q_1\left(x_1,z_1\right)}{\der z_1}\frac{\der}{\der
y_1}\right)\circ\varphi^*$. Now, Theorem \ref{thm_LND} implies
that $\delta_0(z_2)$ is divisible by $x_2^{n_2}$. Since
$\delta_0(z_2)=a\alpha^{-1}(x_2-b)^{n_1}$, we must have $b=0$ and
$n_1\geq n_2$. This proves the first part of the corollary.

Moreover, repeating this analysis with $\varphi^{-1}$ instead of
$\varphi$, we also obtain $n_2\geq n_1$ and so $n_1=n_2=n$.

Since $\varphi:X_{Q_1,n}\to X_{Q_2,n}$ is a morphism, we know that
$\varphi^*\left(x^ny-Q_2(x,z)\right)$  belongs to the ideal
$\left(x^ny-Q_2(x,z)\right)$. In particular, when $x=0$, it
implies that $Q_2\left(0,\alpha
z+\beta(0)\right)=\varphi^*\left(Q_2(0,z)\right)\in\left(Q_1(0,z)\right)$.
Thus $\deg(Q_2(0,z))\geq\deg(Q_1(0,z))$.

Working with $\varphi^{-1}$, the same analysis allows us to
conclude that $\deg(Q_1(0,z))\geq\deg(Q_2(0,z))$. Moreover, it
implies that $Q_2\left(0,\alpha z+\beta(0)\right)=\mu Q_1(0,z)$
for a certain constant $\mu\in\C^*$.

Since the case $n_1=n_2=1$ was already done by Daigle
\cite{Daigle}, this suffices to prove the second part of the
corollary.

It remains to prove the third part.

Let $X_{Q_1,n}$ and $X_{Q_2,n}$ be two equivalent Danielewski
hypersurfaces with $n\geq2$, and let $\Phi$ be an algebraic
automorphism of $\C^3$ such that $\Phi(X_{Q_1,n})=X_{Q_2,n}$.

Since the polynomial $x^ny-Q_1(x,z)$ is irreducible, there exists
a nonzero constant $\mu\in\C^*$ so that
$\Phi^*(x^ny-Q_2(x,z))=\mu(x^ny-Q_1(x,z))$.

Thus, $\Phi$ induces an isomorphism $\Phi_c$ between the
Danielewski hypersurfaces of equation $x^ny-Q_2(x,z)=\mu c$ and
$x^ny-Q_1(x,z)=c$ for every $c\in\C$.

Since $n\geq 2$, the Makar-Limanov invariant of these
hypersurfaces is $\C[x]$. By (1), we obtain now that the image by $\Phi^*$ of the ideal
$(x,x^ny-Q_2(x,z)-\mu c)$ belongs to the ideal
$(x,x^ny-Q_1(x,z)-c)=(x,Q_1(0,z)+c)$ for each $c\in\C$. It turns out that
$$\Phi^*(x)\in\bigcap_{c\in\C}(x,Q_1(0,z)+c)=(x).$$
\noindent Since $\Phi$ is invertible, this implies that
$\Phi^*(x)=ax$ for a certain constant $a\in\C^*$. Thus
$$-\mu Q_1(0,z)\equiv \mu(x^ny-Q_1(x,z))\equiv\Phi^*(x^ny-Q_2(x,z))\equiv
-Q_2(0,\Phi^*(z))\mod{(x)}.$$ \noindent Since $\deg Q_1(0,z)=\deg
Q_2(0,z)$ (by the second part of the corollary), this implies that
$\Phi^*(z)\equiv\alpha z+\beta\mod{(x)}$ for certain constants
$\alpha$ and $\beta$ such that $Q_2(0,\alpha z+\beta)=\mu
Q_1(0,z)$.

Thus, we can write $\Phi^*(z)=\alpha z+\beta+xB(x,y,z)$, where $B$
is polynomial of $\C[x,y,z]$. 

Now, we use again the first part of the corollary. For every $c\in\C$,
there exist a constant $\alpha_c\in\C^*$ and a polynomial
$\beta_c\in\C^{[1]}$ such that $$\Phi^*(z)=\alpha
z+\beta+xB(x,y,z)\equiv \alpha_cz+\beta_c(x)
\mod(x^ny-Q_1(x,z)-c).$$ Therefore, for every $c\in\C$, we have
$\alpha_c=\alpha$, $\beta_c(0)=\beta$  and $$B(x,y,z)\equiv x^{-1}(\beta_c(x)-\beta)\mod(x^ny-Q_1(x,z)-c).$$

In particular $B$ has the following property: For infinitely many constants $c\in\C$, there exist polynomials $r_c(x)\in\C[x]$ and $s_c(x,y,z)\in\C[x,y,z]$ such that $$B(x,y,z)=r_c(x)+s_c(x,y,z)(x^ny-Q_1(x,z)-c).$$ We will show that any polynomial with this property must  belong to $\C[x,x^ny-Q_1(x,z)]$. Remark that it suffices to show that at least one polynomial $s_c$  belongs to $\C[x,x^ny-Q_1(x,z)]$.

In order to see this, we define a degree function  $d$ on
$\C[x,y,z]$ by posing, for every
$f\in\C[x,y,z]$, $d(f):=\deg(\tilde
f(y,z))$, with $\tilde
f(y,z)=f(x,y,z)\in\C[x][y,z]$.

Let $B(x,y,z)=r_{c_0}(x)+s_{c_0}(x,y,z)(x^ny-Q_1(x,z)-c_0)$ for one $c_0\in\C$.  Then, $s_{c_0}$ satisfies also the above property and its degree $d(s_0)$ is strictly less than $d(B)$.

Therefore, the desired result can be obtained by decreased induction on the degree $d$.
\end{proof}

\section{Equivalence classes}

In this section, we prove the following result.

\begin{thm}\label{thm:classification-plgts}
Two Danielewski hypersurfaces $X_{Q_1,n_1}$ and $X_{Q_2,n_2}$ are
equivalent if and only if  $n_1=n_2=n$ and there exist  $a,\alpha,
\mu\in\C^*$, $\beta\in\C$ and $B\in\C^{[2]}$ such that
$$Q_2(ax,\alpha z+ \beta+xB(x,Q_1(x,z)))\equiv\mu
Q_1(x,z)\mod{(x^n)}.$$
\end{thm}

\begin{rem}
We will show in the next section (Proposition
\ref{prop:plgts_non_equiv}) that this theorem implies that every
Danielewski hypersurface $X_{Q,n}$ with $n\geq2$ admits at least
two non-equivalent embeddings into $\C^3$.
\end{rem}

Before proving Theorem \ref{thm:classification-plgts}, let us give
another result. Given two Danielewski hypersurfaces,
it is not easy to check if  the second condition in Theorem
\ref{thm:classification-plgts} is fulfilled. Therefore, we  also
show that any Danielewski hypersurface is equivalent to another
one which is unique up to an affine automorphism.

\begin{thm}\label{thm:plgt-normal}\hfill
\begin{enumerate}
\item Every Danielewski hypersurface is equivalent to a
Danielewski hypersurface $X(p,\{q_i\}_{i=2..\deg(p)},n)$ defined
by an equation of the form
$$x^ny-p(z)-x\sum_{i=2}^{\deg(p)}p^{(i)}(z)q_{i}(x,p(z))\quad\text{with}\quad\deg_x(q_i)<n-1.$$
Moreover, there is an algorithmic procedure which computes, given
a Danielewski hypersurface $X$,  a hypersurface
$X(p,\{q_i\}_{i=2..\deg(p)},n)$ which is equivalent to $X$.

\item Two such Danielewski hypersurfaces
$X(p_1,\{q_{1,i}\}_{i=2..\deg(p_1)},n_1)$ and
$X(p_2,\{q_{2,i}\}_{i=2..\deg(p_2)},n_2)$ are equivalent if and
only if $n_1=n_2$, $\deg(p_1)=\deg(p_2)=d$ and there exist some
constants $a,\alpha,\mu\in\C^*$, $\beta\in\C$ such that
$p_1(\alpha z+\beta)=\mu p_2(z)$ and $a\alpha^{-i}q_{2,i}(ax,\mu
t)=q_{1,i}(x,t)$ for every $2\leq i\leq d$.
\end{enumerate}
\end{thm}

\begin{rem}
This result generalizes the classification of Danielewski
hypersurfaces of the form $x^2y-z^2-xq(z)$  given by
Moser-Jauslin and the author in \cite{MJ-P}.
\end{rem}

\begin{proof}[Proof of Theorem \ref{thm:classification-plgts}]
Let $X_{Q_1,n_1}$ and $X_{Q_2,n_2}$ be two equivalent Danielewski
hypersurfaces. Then, the second part of Corollary
\ref{cor:LND} implies  $n_1=n_2=n$.

If $n=1$, the result is already known. Indeed, by Lemma \ref{lem:mod_x^n}, every Danielewski hypersurface $X_{Q,1}$ with
$n=1$ is equivalent to a one of the form $X_{p,1}$ with
$p(x,z)=p(z)\in\C[z]$. Then, Daigle \cite{Daigle} has proven that two
such hypersurfaces $X_{p_1,1}$ and $X_{p_2,1}$ are isomorphic if
and only if $p_2(az+b)=\mu p_1(z)$ for some constants
$a,\mu\in\C^*$ and $b\in\C$.

Now, assume $n\geq2$ and let $\Phi$ be an  automorphism of $\C^3$
such that $\Phi(X_{Q_1,n})=X_{Q_2,n}$.  Corollary \ref{cor:LND},
gives us constants $a,\alpha\in\C^*$, $\beta\in\C$ and a
polynomial $B\in\C^{[2]}$ such that $\Phi^*(x)=ax$ and
$\Phi^*(z)=\alpha z+\beta+xB(x,x^ny-Q_1(x,z))$. Since the
polynomial $x^ny-Q_1(x,z)$ is irreducible, there exists a nonzero
constant $\mu\in\C^*$ so that
$\Phi^*(x^ny-Q_2(x,z))=\mu(x^ny-Q_1(x,z))$. It turns out that
$Q_2(ax,\alpha z+ \beta+xB(x,-Q_1(x,z)))\equiv\mu
Q_1(x,z)\mod{(x^n)}$, as desired.

Conversely, let $X_{Q_1,n}$ and $X_{Q_2,n}$ be two Danielewski
hypersurfaces with  $Q_2(ax,\alpha z+
\beta+xB(x,Q_1(x,z)))\equiv\mu Q_1(x,z)\mod{(x^n)}$ for some
$a,\alpha, \mu\in\C^*$, $\beta\in\C$ and $B\in\C^{[2]}$.

We pose $$R(x,y,z)=x^{-n}\left(Q_2\left(ax,\alpha z+
\beta+xB\left(x,-x^ny+Q_1\left(x,z\right)\right)\right)-\mu
Q_1\left(x,z\right)\right)\in\C[x,y,z]$$ \noindent and define an
endomorphism of $\C^3$ by
$$\Phi(x,y,z)=(ax,a^{-n}\mu y+a^{-n}R(x,y,z),\alpha z+
\beta+xB(x,-x^ny+Q_1(x,z))).$$

Remark that $\Phi^*(x^ny-Q_2(x,y))=\mu (x^ny-Q_1(x,y))$.
Therefore, the theorem will be proved if we show that $\Phi$ is
invertible.

It suffices to prove that $\Phi^*$ is surjective, i.e.
$$\C[x,y,z]\subset\Phi^*(\C[x,y,z])=\C[\Phi^*(x),\Phi^*(y),\Phi^*(z)].$$

We know already  that $x$ and $P_1:=x^ny-Q_1(x,z)$ are in the image
of $\Phi^*$.

Then, since $z=\alpha^{-1}(\Phi^*(z)-\beta-xB(x,-P_1))$, we obtain
that $z$ belongs to $\Phi^*(\C[x,y,z])$.

It remains to show that $y$ belongs to the image of $\Phi^*$. To
do this, we first remark that there exist  polynomials
$f,g\in\C^{[3]}$ such that
$$R(x,y,z)=xyf(x,z,P_1)+g(x,z,P_1).$$

Thus, $y(\mu+xf(x,z,P_1))=a^n\Phi^*(y)-g(x,z,P_1)$
$\in\Phi^*(\C[x,y,z])$. Now, choose some polynomials  $\tilde f$
and $\tilde g$ such that
$$(\mu+xf(x,z,P_1))\tilde f(x,z,P_1)=1+x^n\tilde
 g(x,z,P_1).$$ Then, we can write
\begin{align*}
y(\mu+xf(x,z,P_1))\tilde  f(x,z,P_1)&=y(1+x^n\tilde
 g(x,z,P_1))\\ &=y+(P_1-Q_1(x,z))\tilde g(x,z,P_1).
 \end{align*}
This implies that $y$ belongs to the image of $\Phi^*$ and proves the theorem.
\end{proof}

\begin{proof}[Proof of Theorem \ref{thm:plgt-normal}]
Let $X_{Q,n}$ be a Danielewski hypersurface. Let $p$ and $q$
denote the polynomials such that $Q(x,z)=p(z)+xq(x,z)$.

We can write $q(x,z)$ in the following form:
$q(x,z)=\sum_{i=1}^{\deg(p)}p^{(i)}(z)q_i(x,p(z))$ for some
polynomials $q_1,\ldots,q_{\deg(p)}\in\C^{[2]}$.

By Lemma \ref{lem:mod_x^n},  we can assume   $\deg_x(q_i(x,p(z)))<n-1$
for every index $1\leq i\leq\deg(p)$.

Now, rewrite
$$Q(x,z)=p(z)+\sum_{k=1}^{n-1}x^k\sum_{i=1}^{\deg(p)}p^{(i)}(z)q_{i,k}(p(z))$$
\noindent  for suitable polynomials $q_{i,k}\in\C^{[1]}$. Let
$1\leq k_0\leq n-1$ be a fixed integer.  Then,
$$Q\left(x,z-x^{k_0}q_{1,k_0}\left(p\left(z\right)\right)\right)\equiv  p(z)+\sum_{k=1}^{k_0}x^k\sum_{i=1}^{\deg(p)}p^{(i)}(z)q_{i,k}(p(z))-x^{k_0}q_{1,k_0}(p(z))p'(z)\mod(x^{k_0+1}).$$
Therefore, we obtain, using Theorem
\ref{thm:classification-plgts}, that $X_{Q,n}$ is equivalent to a
hypersurface of equation
$$x^ny=p(z)+\sum_{k=1}^{n-1}x^k\sum_{i=1}^{\deg(p)}p^{(i)}(z)\tilde q_{i,k}(p(z))$$
\noindent with $\tilde q_{1,k_0}=0$ and $\tilde q_{i,k}=q_{i,k}$
if $(i,k)\in[1,k_0]\times[1,\deg(p)]\setminus\{(1,k_0)\}$. Then, it
is easy to prove by induction on $k_0$ that $X_{Q,n}$ is
equivalent to a hypersurface of the desired form.

Now, we will prove the second part of the theorem.

Let $X_j=X(p_j,\{q_{j,i}\}_{i=2..\deg(p_j)},n_j)$, $j=1,2$,  and pose
$Q_j=p_j(z)+x\sum_{i=2}^{\deg(p_j)}p_j^{(i)}(z)q_{j,i}(x,p_j(z))$.

If $X_1$ and $X_2$ are equivalent, then, it follows from
Theorem \ref{thm:classification-plgts}, that $n_1=n_2=n$ and that
there exist  $a,\alpha, \mu\in\C^*$, $\beta\in\C$ and
$B\in\C^{[2]}$ such that
$$Q_2(ax,\alpha z+ \beta+xB(x,Q_1(x,z)))\equiv\mu
Q_1(x,z)\mod{(x^n)}.$$

This implies  $p_2(\alpha z+\beta)=\mu p_1(z)$. Thus
$\deg(p_1)=\deg(p_2)=d$.

First, we prove that $B(x,\cdot)\equiv 0\mod(x^{n-1})$. In
order to do this, suppose that we can write $B(x,t)\equiv
b_k(t)x^k$ for some $0\leq k\leq n-2$ and
$b_k(t)\in\C[t]\setminus{0}$. Then, we obtain the following
equalities modulo $(x^{k+2})$.

\begin{align*}
\mu Q_1(x,z)&\equiv Q_2\left(ax,\alpha
z+\beta+xB\left(x,Q_1\left(x,z\right)\right)\right) \mod(x^{k+2})
\\
&\equiv Q_2\left(ax,\alpha
z+\beta+x^{k+1}b_k\left(p_1\left(z\right)\right)\right)  \\
&\equiv p_2\left(ax,\alpha
z+\beta+x^{k+1}b_k\left(p_1\left(z\right)\right)\right)+ax\sum_{i=2}^dp_2^{(i)}(\alpha z+\beta)q_{2,i}(ax,p_2(\alpha z+\beta))\\
&\equiv \mu p_1(z)+x^{k+1}b_k(p_1(z))(p_2)'(\alpha z+\beta)+ax\sum_{i=2}^dp_2^{(i)}(\alpha z+\beta)q_{2,i}(ax,\mu p_1(z))\\
&\equiv \mu p_1(z)+x^{k+1}b_k(p_1(z))\alpha^{-1}\mu
p_1'(z)+ax\sum_{i=2}^d\alpha^{-i}\mu p_1^{(i)}(z)q_{2,i}(ax,\mu
p_1(z)).
\end{align*}

This would imply
$$\alpha^{-1}(p_1)'(z)b_k(p_1(z))\equiv\frac{\sum_{i=2}^{d}\left(p_1^{(i)}(z)q_{1,i}(x,p_1(z))-a\alpha^{-i} p_1^{(i)}(z)q_{2,i}(ax,\mu
p_1(z)) \right)}{x^k}\mod(x^{k+1}),$$ \noindent what is impossible.

Therefore, $B(x,\cdot)\equiv 0\mod(x^{n-1})$. Since, by
hypothesis,  $\deg_x(Q_j(x,z))<n$ for $j=1,2$, it follows that
$Q_2(ax,\alpha z+ \beta)=\mu Q_1(x,z)$. Then, we can easily check
that this last equality implies $a\alpha^{-i}q_{2,i}(ax,\mu
t)=q_{1,i}(x,t)$ for every $2\leq i\leq d$, as desired. This
concludes the proof.
\end{proof}

\section{Standard forms}

In \cite{Dubouloz-Poloni}, A. Dubouloz and the author proved that
every   Danielewski hypersurface $X_{Q,n}$  where  $Q(x,z)$ is a polynomial
such that $Q(0,z)$ has $d\geq2$ simple roots,  is isomorphic to a
hypersurface of a certain type (called standard form) and then
classified all these standard forms up to isomorphism.

In this section, we will generalize these results even when $Q(0,z)$ has multiple roots. In order to do this, we first generalize the definition of standard form  given in \cite{Dubouloz-Poloni}.

\begin{defi}
We say that a Danielewski hypersurface $X_{Q,n}$ is in
\emph{standard form} if the polynomial $Q$ can be written as
follows:
$$Q(x,z)=p(z)+xq(x,z),\quad\text{with }\deg_z(q(x,z))<\deg(p).$$
\end{defi}

We also introduce a notion of reduced standard form.

\begin{defi}
A Danielewski hypersurface $X_{Q,n}$ is in \emph{reduced standard
form} if $\deg_x\left(Q\left(x,z\right)\right)<n$ and
$$Q(x,z)=p(z)+xq(x,z),\quad\text{with }\deg_z(q(x,z))<\deg(p)-1.$$
\end{defi}

When $X_1$ and $X_2$ are two isomorphic Danielewski hypersurfaces
with $X_2$ in (reduced) standard form, we say that $X_2$ is a
(reduced) standard form of $X_1$.

\begin{example}\hfill
\begin{enumerate}
 \item  Danielewski hypersurfaces defined by equations of the form $x^ny-p(z)=0$ are in reduced standard form (These hypersurfaces were studied by  Makar-Limanov in \cite{ML01});
 \item The Danielewski hypersurfaces,  studied by Danielewski \cite{Danielewski} and Wilkens \cite{Wilkens}, defined by $x^2y-z^2-h(x)z=0$ are in standard form;
\item  Danielewski hypersurfaces $X_{\sigma,n}$ defined by
equations $x^ny=\prod_{i=1}^d(z-\sigma_i(x))$, where
$\sigma=\{\sigma_i(x)\}_{i=1\cdots d}$ is a collection of $d\geq2$
polynomials, are in standard form; (They are those we
have called in standard form in \cite{Dubouloz-Poloni})\item If
$r(x)\in\C[x]$ is a non constant polynomial, then a Danielewski
hypersurface defined by  $x^ny-r(x)p(z)=0$ is not in standard
form. (They were studied by  Freudenburg and Moser-Jauslin
\cite{Freudenburg-MJ}.)
\end{enumerate}
\end{example}

We will now prove that every Danielewski hypersurface is
isomorphic to a one in reduced standard form. Our proof will be
based on the following lemma which comes from \cite{Freudenburg-MJ}.

\begin{rem}
Since every Danielewski hypersurface is isomorphic to a one in reduced standard form, the notion of reduced standard form is, in some sense,  more relevant than the notion of standard form if we are interested in the classification of  Danielewski hypersurfaces. Nevertheless, the notion of standard form has an interest too. Indeed, nice properties are true for every Danielewski hypersurfaces in (not necessarily reduced) standard forms. For example, we will see that all their automorphisms extend to automorphisms of the ambient space. Recall that this does not hold for all Danielewski hypersurfaces. (see \cite{Dubouloz-Poloni})
\end{rem}

\begin{lem}\label{lem:Lucy-Gene}
Let  $n\geq1$ be a  natural number and $Q_1(x,z)$ and $Q_2(x,z)$
be two polynomials of $\C[x,z]$ such that
$$Q_2(x,z)=(1+x\pi(x,z))Q_1(x,z)+x^{n}R(x,z)$$

\noindent for some polynomials $\pi(x,z),R(x,z)\in\C[x,z]$.

Then, the endomorphism of $\C^3$ defined by
$$\Phi(x,y,z)=\left(x , \left(1+x\pi\left(x,z\right)\right)y+R\left(x,z\right) , z\right)$$
\noindent induces an isomorphism  $\varphi:X_{Q_1,n}\to
X_{Q_2,n}$.
\end{lem}

\begin{proof}
Remark that, since
$\Phi^*\left(x^ny-Q_2\left(x,z\right)\right)=\left(1+x\pi\left(x,z\right)\right)\left(x^ny-Q_1\left(x,z\right)\right)$,
\noindent $\Phi$ induces a morphism $\varphi:X_{Q_1,n}\to
X_{Q_2,n}$.

Let  $f(x,z)$ and $g(x,z)$ be two polynomials in $\C[x,z]$ so that
$(1+x\pi(x,z))f(x,z)+x^ng(x,z)=1$ and define $\Psi$, an
endomorphism of $\C^3$, by posing
$$\left\{ \begin{array}{l}
\Psi^*(x)=x \\
\Psi^*(y)=f(x,z)y+g(x,z)Q_1(x,z)-f(x,z)R(x,z) \\
\Psi^*(z)=z
\end{array} \right. $$
We check easily that
$$\Psi^*\left(x^ny-Q_1\left(x,z\right)\right)=f\left(x,z\right)\left(x^ny-Q_2\left(x,z\right)\right)$$
\noindent and that
$$\begin{array}{l}
\Psi^*\circ\Phi^*(x)=x \ ; \\
\Psi^*\circ\Phi^*(y)=y-g(x,y)\left(x^ny-Q_2\left(x,z\right)\right) \ ; \\
\Psi^*\circ\Phi^*(z)=z. \\
\end{array}$$

Therefore, the  restriction of $\Psi^*\circ\Phi^*$ to
$S_{Q_2,n}=\C[X_{Q_2,n}]$ is identity. Hence,  $\Psi$ induces the inverse
morphism of $\varphi$, and $X_{Q_1,n}\simeq X_{Q_2,n}$.
\end{proof}

\begin{thm}\label{thm:formesstandards}
Every Danielewski hypersurface is isomorphic to a Danielewski
hypersurface in reduced standard form. Furthermore, there is an
algorithmic procedure which computes one of the reduced standard
forms of a given Danielewski hypersurface.
\end{thm}

\begin{proof}
Let $X=X_{Q,n}$ be a Danielewski hypersurface and denote
$Q(x,z)=p(z)+xq(x,z)$ with $p(z)\in\C[z]$ and $q(x,z)\in\C[x,z]$.

One can construct, by induction on $m\geq0$, two polynomials
 $q_{s,m}(x,z)$ and
$\pi_m(x,z)$ so that $\deg_z(q_s(x,z))<\deg(p)$  and $Q(x,z)
\equiv (1+x\pi_{m}(x,z))(p(z)+xq_{s,m}(x,z)) \mod{(x^{m+1})}$.

\noindent Indeed, this assertion is obvious for $m=0$, whereas, if
it is true for a rank $m$, we can write:

\begin{align*}
p(z)+xq(x,z)& \equiv
(1+x\pi_{m}(x,z))(p(z)+xq_{s,m}(x,z)) \mod{(x^{m+1})}\\
 & =(1+x\pi_{m}(x,z))(p(z)+xq_{s,m}(x,z))+x^{m+1}R_{m+1}(x,z) \\
& \equiv (1+x\pi_{m}(x,z)+x^m\tilde{\pi}_{m+1}(z))(p(z)+xq_{s,m}(x,z)+x^{m}r_{m+1}(z)) \mod{(x^{m+1})} \\
& \equiv (1+x\pi_{m+1}(x,z))(p(z)+xq_{s,m+1}(x,z)) \mod{(x^{m+1})},\\
 \end{align*}

\noindent where
$R_{m+1}(0,z)=p(z)\tilde{\pi}_{m+1}(z)+{r}_{m+1}(z)$ is the
Euclidean division (in  $\C[z]$) of $R_{m+1}(0,z)$ by $p$.

Thus, we obtain $$p(z)+xq(x,z)=
(1+x\pi_{n-1}(x,z))(p(z)+xq_{s,n-1}(x,z))+x^nR_n(x,z).$$  Lemma
\ref{lem:Lucy-Gene} allows us to conclude that $X$ is isomorphic
to the Danielewski hypersurface in standard form $X_s$ defined by
the equation $x^ny-p(z)-xq_{s,n-1}(x,z)=0$.

In order to obtain a reduced standard form, we rewrite
$$p(z)+xq_{s,n-1}(x,z)=\sum_{i=0}^{d}a_iz^i+x\sum_{i=0}^{d-1}z^i\alpha_i(x)$$
\noindent and consider the automorphism of $\C^3$ defined by
$$\Phi:(x,y,z)\mapsto(x,y,z-x(da_d)^{-1}\alpha_{d-1}(x)).$$

One checks that the polynomial $\Phi^*(x^ny-p(z)-xq_{s,n-1}(x,z))$
satisfies the second condition in the definition of Danielewski
hypersurface in reduced standard form. Finally, the first condition
can be obtain easily by applying Lemma \ref{lem:Lucy-Gene}.

This proof gives an algorithm for finding a (reduced) standard
form of given Danielewski hypersurface.
\end{proof}

It should be noticed that a Danielewski hypersurface is in general not equivalent to its (reduced) standard form given by Theorem \ref{thm:formesstandards}. Morover, one can use this fact to construct  non-equivalent embeddings for every Danielewski hypersurface of non-trivial Makar-Limanov invariant.

\begin{prop}\label{prop:plgts_non_equiv}
Every Danielewski hypersurface $X_{Q,n}$ with $n\geq2$ admits at
 least two non-equivalent embeddings into $\C^3$.
\end{prop}

\begin{proof}
Since, by Theorem \ref{thm:formesstandards}, every Danielewski
hypersurface is isomorphic to a one in standard form, it suffices
to show that every Danielewski hypersurface in standard form
$X_{Q,n}$ with $n\geq2$ admits at
 least two non-equivalent embeddings in $\C^3$.

 Let $X=X_{Q,n}$ be a Danielewski hypersurface in standard form with
$n\geq2$. Then, due to Lemma \ref{lem:Lucy-Gene}, $X$ is
isomorphic to the hypersurface $Y=X_{(1+x)Q(x,z),n}$.
Nevertheless, it turns out that $X$ and $Y$ are non-equivalent
hypersurfaces of $\C^3$. Indeed, if they were,  Theorem
\ref{thm:classification-plgts} would give us constants $a,\alpha,
\mu\in\C^*$, $\beta\in\C$ and a polynomial $B\in\C^{[2]}$ such that
\begin{equation*}\label{equ}
    (1-ax)Q(ax,\alpha z+ \beta+xB(x,Q(x,z)))\equiv\mu
Q(x,z)\mod{(x^n)}.
\end{equation*}

In turn, if we denote $Q(x,z)=p(z)+xq(x,z)$, it would lead the following equalities modulo $(x^2)$:
\begin{align*} \mu Q(x,z) &\equiv \mu
(p(z)+xq(0,z))\\&\equiv(1-ax)Q(ax,\alpha z+
\beta+xB(0,Q(0,z)))
\\&\equiv(1-ax)(p\big(\alpha z+\beta+xB(0,p(z)))+xq(0,\alpha
z+\beta))
\\&\equiv p(\alpha z+\beta)+x(B(0,p(z))p'(\alpha
z+\beta)+q(0,\alpha z+\beta)-ap(\alpha z+\beta)).\\
\end{align*}
Thus $$B(0,p(z))p'(\alpha z+\beta)+q(0,\alpha
z+\beta)-ap(\alpha z+\beta)=\mu q(0,z)$$\noindent which is impossible since
$\deg(q(0,z))<\deg(p)$ by definition of a standard form.
\end{proof}

\begin{rem}
This proof is similar to
 the proof of Freudenburg and Moser-Jauslin in \cite{Freudenburg-MJ} for hypersurfaces of equation $x^ny=p(z)$ with $n\geq2$. In their article, they also have constructed non-equivalent embeddings into $\C^3$ for Danielewski hypersurfaces of
 the form
$xy-z^{d}-1=0$ for some $d\in\N$. Nevertheless,  we do not know if every
Danielewski hypersurface $X_{Q,1}$ admits non equivalent embeddings
into $\C^3$. For instance,  the following question , which they
posed in \cite{Freudenburg-MJ}, is still open.
\end{rem}

\begin{question}
Does the hypersurface of equation $xy+z^2=0$ admit a unique embedding into $\C^3$?
\end{question}

Remark also that the two non-equivalent embeddings of a
Danielewski hypersurface $X_{Q,n}$ with $n\geq2$ which we
construct
 in Proposition
\ref{prop:plgts_non_equiv}  are \emph{analytically} equivalent.
 Indeed,  it
can be  easily seen, as  in \cite{Freudenburg-MJ} and \cite{Dubouloz-Poloni}, that a Danielewski hypersurface is
analytically equivalent to its standard form given by
Theorem \ref{thm:formesstandards}. In turn, we obtain the
following result.

\begin{prop}
If $X_1$ and $X_2$ are two isomorphic  Danielewski hypersurfaces,
then there is an analytic automorphism $\Psi$ of $\C^3$ such that
$\Psi(X_1)=X_2$.
\end{prop}

\begin{proof}
Let $X=X_{Q,n}$ be a Danielewski hypersurface and let $X_{Q_s,n}$
be its standard form given by the theorem
\ref{thm:formesstandards}. In the proof of this theorem, we have
seen that $Q(x,z)=(1+x\pi(x,z))Q_s(x,z)+x^nR(x,z)$ for certain
polynomials $\pi(x,z),R(x,z)\in\C[x,z]$. Now, consider the
following analytic automorphism of $\C^3$.
$$\Psi:(x,y,z)\mapsto(x,\exp(xf(x,z))y-x^{-n}(\exp(xf(x,z))-1-x\pi(x,z))Q_s(x,z)+R(x,z),z),$$
\noindent where $f(x,z)\in\C[x,z]$ is a polynomial so that
$\exp(xf(x,z))\equiv1+x\pi(x,z)\mod(x^n)$. One checks that
$\Psi^*(x^ny-Q(x,z))=x^ny-Q_s(x,z)$. Thus, $\Psi$ maps $X_{Q,n}$
to its standard form $X_{Q_s,n}$. Then, the result follows from
Proposition \ref{prop:ext-isos}, which will be proved at the end
of this paper.
\end{proof}

\section{Classification up to isomorphism}

Finally, we  give the classification of Danielewski hypersurfaces
in standard form. Together with the theorem
\ref{thm:formesstandards}, this effectively classifies all the
Danielewski hypersurfaces up to isomorphism of algebraic
varieties.

\begin{thm}\label{thm:class-standard}\hfill
\begin{enumerate}
\item Two Danielewski hypersurfaces $X_{Q_1,n_1}$ and
$X_{Q_2,n_2}$ in standard form are isomorphic if and only if the
two following conditions are satisfied:
\begin{enumerate}
\item $n_1=n_2=n$; \item $\exists a,\alpha,
\mu\in\C^*\quad\exists\beta(x)\in\C[x]$ such that $Q_2(ax,\alpha z
+\beta(x))\equiv\mu Q_1(x,z) \mod(x^n)$.
\end{enumerate}
\item Two Danielewski hypersurfaces $X_{Q_1,n_1}$ and
$X_{Q_2,n_2}$ in reduced standard form are isomorphic if and only
if the two following conditions are satisfied:
\begin{enumerate}
\item $n_1=n_2$; \item $\exists a,\alpha,
\mu\in\C^*\quad\exists\beta\in\C$ such that $Q_2(ax,\alpha z
+\beta)=\mu Q_1(x,z)$.
\end{enumerate}
\end{enumerate}
\end{thm}

\begin{proof}
 Let $X_1=X_{Q_1,n_1}$ and $X_2=X_{Q_2,n_2}$ be two
isomorphic Danielewski hypersurfaces in standard form and let
$\varphi:X_1\to X_2$ be an isomorphism. Then Corollary
\ref{cor:LND} implies that $n_1=n_2=n$. Since the case $n=1$ was
already done by Daigle \cite{Daigle}, we can suppose that
$n\geq2$.

Denote by $x_i,y_i,z_i$ the images of $x,y,z$ in the coordinate
ring $\C[X_i]$ for $i=1,2$.  Then, due to Corollary \ref{cor:LND},
there exist constants $a,\alpha\in\C^*$ and a polynomial
$\beta(x)\in\C[x]$ such that $\varphi^*(x_2)=ax_1$ and
$\varphi^*(z_2)=\alpha z_1+\beta(x_1)$.

Moreover,  we have proven in the proof of Corollary
\ref{cor:LND}, that $Q_2\left(0,\alpha z_2+\beta(0)\right)=\mu
Q_1(0,z_1)$ for a certain constant $\mu\in\C^*$.

 Thus, viewing $\C[X_i]$ as a subalgebra of
$\C[x_i,x_i^{-1},z_i]$ with $y_i=x_i^{-n}Q_i(x_i,z_i)$, we obtain
$$\varphi^*(y_2)=\varphi^*(x_2^{-n}Q_2(x_2,z_2))=(ax_1)^{-n}Q_2(ax_1,\alpha
z_1+\beta(x_1))=\mu a^{-n}y_1+(ax_1)^{-n}\Delta(x_1,z_1),$$
\noindent where $\Delta(x_1,z_1)=Q_2(ax_1,\alpha
z_1+\beta(x_1))-\mu Q_1(x_1,z_1)$.

Remark that $\deg_{z_1}\Delta(x_1,z_1)<\deg_{z_1}Q_1(0,z_1)$ since
$X_1$ and $X_2$ are in standard form.

It turns out that $x_1^{-n}\Delta(x_1,z_1)\in\C[x_1,z_1]$ since
any polynomial of $\C[X_1]\subset\C[x_1^\pm,z_1]$ with negative
degree in $x_1$ has obviously a degree in $z_1$ at least equal to
$\deg_{z_1}Q_1(0,z_1)$. Thus, $\Delta(x,z)\equiv0\mod(x^n)$ and
$X_1$ and $X_2$ fulfill conditions (1) (a) and (1) (b).

If $X_1$ and $X_2$ are in reduced standard form, then we see easily that
$\Delta(x,z)\equiv0\mod(x^n)$ is possible only if  $\beta(x)\equiv
\beta(0)\mod(x^n)$. If so $Q_2(ax_1,\alpha z_1+\beta(0))=\mu
Q_1(x_1,z_1)$ and $X_1$ and $X_2$ fulfill the conditions (1) (a)
and (2) (b).

Conversely, suppose that $X_1=X_{Q_1,n}$ $X_2=X_{Q_2,n}$ are two
Danielewski hypersurfaces which satisfy the conditions (a) and (b)
of part (1). Then the following triangular automorphism of $\C^3$
induces  an isomorphism between $X_1$ and $X_2$:
$$(x,y,z)\mapsto(ax,\mu a^{-n}y+(ax)^{-n}(Q_2(ax,\alpha z+\beta(x))-\mu
Q_1(x,z)),\alpha z+\beta).$$
\end{proof}

As a corollary, we observe that two isomorphic Danielewski
hypersurfaces in standard form are equivalent via a triangular
automorphism of $\C^3$, and that two isomorphic Danielewski
hypersurfaces in reduced standard form are equivalent via an
affine one. In fact, we have even proven a stronger result in the
proof of Theorem \ref{thm:class-standard}.

\begin{prop}\label{prop:ext-isos}
Every isomorphism between two isomorphic Danielewski hypersurfaces
in standard form can be lifted to a triangular automorphism of
$\C^3$.
\end{prop}

\bibliographystyle{plain}
\bibliography{DanielewskiHypersurfacesReferences}

\end{document}